\documentclass{amsproc}
\usepackage{graphicx}
\usepackage{amscd}
\usepackage{amsmath}
\usepackage{amsfonts}
\usepackage{amssymb}
\theoremstyle{plain}

\newtheorem{definition}{Definition}

\newtheorem{lemma}{Lemma}

\newtheorem{proposition}{Proposition}
\newtheorem{remark}{Remark}

\newtheorem{theorem}{Theorem}
\numberwithin{equation}{section}

\begin{document}
\title[Almost Split Sequences]{Local Theory of Almost Split Sequences for Comodules}
\author{William Chin}
\address{DePaul University, Chicago, Illinois 60614}
\email{wchin@condor.depaul..edu}
\urladdr{http://www.depaul.edu/\symbol{126}wchin}
\author{Mark Kleiner}
\address{Syracuse University, Syracuse, NY 13244}
\email{mkleiner@mail.syr.edu}
\urladdr{}
\author{Declan Quinn}
\address{Syracuse University, Syracuse, NY 13244}
\email{dpquinn@mail.syr.edu}
\urladdr{}
\date{}
\subjclass{}
\keywords{}
\dedicatory{}

\begin{abstract}
We show that almost split sequences in the category of comodules over a
coalgebra $\Gamma$ with finite-dimensional right-hand term are direct limits
of almost split sequences over finite dimensional subcoalgebras. \ In previous
work we showed that such almost split sequences exist if the right hand term
has a quasifinitely copresented linear dual. Conversely, taking limits of
almost split sequences over finte-dimensional comodule categories, we then
show that, for countable-dimensional coalgebras, certain exact sequences exist
which satisfy a condition weaker than being almost split, which we call
``finitely almost split''. Under additional assumptions, these sequences are
shown to be almost split in the appropriate category.
\end{abstract}\maketitle

\section{Introduction}

In [CKQ] we approached the problem of proving the existence of almost split
sequences for comodules following the approach used in [ARS]. In this paper we
focus on the local theory and we study how almost split sequences relate to
almost split sequences over finite dimensional subcoalgebras. Rather than the
more functorial approach of [CKQ] and [Tak], we lean here toward using limits
to obtain results from ones known in the finite-dimensional case. The results
make use of techniques using duality and idempotents in the dual algebra. This
approach is amply represented in other recent work on coalgebras, e.g. [DT, Si].

A fundamental, but surprisingly difficult, result in [CKQ] is the fact that
the functor $\ast$ is a duality on quasifinite injectives. We begin by giving
an easy proof of this result using the characterization of the cohom functor
as the direct limit of duals of $\operatorname*{Hom}$'s.

We show in 3.2 that the almost split sequences with finite-dimensional
right-hand terms are direct limits of almost split sequences over some
finite-dimensional subcoalgebras. This follows from our result in 3.1 which
shows that the transpose of a comodule can be expressed as a direct limit of
finite-dimensional transposes.\ 

In section 4, we attempt to construct almost split sequences over a coalgebra
$\Gamma$ from such sequences over finite-dimensional subcoalgebras. We show
that (4.3 Theorem 2), for countable-dimensional coalgebras, starting with a
finite-dimensional comodule on the right, certain exact sequences exist which
satisfy a condition similar to, but weaker than, being almost split. This type
of sequence, which we call ``right finitely almost split'', was investigated
by E. Green and E. Marcos [GM] in a different setting.\ The assumption on
dimension seems to be a mild one, as any indecomposable coalgebra with a
locally countable Ext-quiver (i.e., $\dim_{\mathrm{k}}$Ext$_{-\Gamma}%
^{1}(S,T)$ is countable for all simple comodules $S,T$) has countable
dimension. The dimension assumption arises mainly from the need to construct a
direct system of almost split sequences; to do this (in the proofs of Theorems
2 and 3) we need to work with a chain of finite-dimensional subcoalgebras
whose union is $\Gamma$, a condition equivalent to $\Gamma$ being of countable
dimension. Under additional assumptions these sequences are shown to be almost
split in the quasifinite comodule category in Theorem 2. This result is
predicted by [CKQ, Corollary 4.3] along with our results in Section 1. In
addition, the result from [CKQ] says that the sequences so obtained are almost
split in the comodule category. \ 

The situation in Theorem 2 dualizes. We start with a finite-dimensional
comodule on the left and construct a left finitely almost split sequence in
the category of prorational modules, which consists of inverse limits of
finite-dimensional comodules. These modules form a dual category to the
category of comodules and are developed in 4.2. As for the right-hand variant,
these sequence are almost split under appropriate hypotheses. This result can
predicted from the result (4.3, Theorem 4), that dualizes [CKQ, Corollary 4.3]
by working in the prorational category.

The Mittag-Leffler condition concerns the exactness of the inverse limit. It
is the main tool for showing that the sequences in Section 4 might be almost
split. The relevant special case is discussed briefly in 4.1.

We close by presenting some examples of almost split sequences over path
coalgebras. Our examples include almost split sequences, starting with
finite-dimensional comodules on either the left or right, that may not be in
the finite-dimensional comodule category. One of the examples exhibits how the
direct limit of almost split sequences may fail to be almost split. More
examples of almost split sequences and AR quivers for comodules appear in the
recent article [KS].

\bigskip

\bigskip\ \textbf{Notation }Let $\Gamma$ denote a coalgebra over the fixed
base field k. Set the following

$\mathcal{M}^{\Gamma}$ the category of right $\Gamma$-comodules.

$\mathcal{M}_{f}^{^{\Gamma}}$ the category of finite-dimensional right
$\Gamma$ -comodules

$\mathcal{M}_{q}^{^{\Gamma}}$ the category of quasifinite right $\Gamma$ -comodules

$\mathcal{M}_{qc}^{^{\Gamma}}$ the category of quasifinitely
copresented\ right $\Gamma$ -comodules

$\mathcal{I}^{\Gamma}$ the category of quasifinite injective left $\Gamma$-comodules

$R=$\textrm{D}$\Gamma$\ the dual algebra $Hom_{k}(\Gamma,$k$)$

$\mathcal{M}_{R}$ the category of right $R$-modules

$\mathcal{D}^{\Gamma}$ the category of duals of left $\Gamma$ -comodules

$\mathcal{D}_{q}^{\Gamma}$ the category of duals of quasifinite left $\Gamma$ -comodules

\bigskip

$h_{-\Gamma}(\_,\_)$ the cohom functor

$\square$ the cotensor product (over $\Gamma$)

D the linear dual $\operatorname*{Hom}_{k}(\_,$k$)$

$(-)^{\ast}$ the functor $h_{-\Gamma}(\_,\Gamma).$

\bigskip

We will use the obvious left-handed variants of these notations$.$ We shall
generally follow the conventions of [CKQ] an the books [Mo, Sw]. The socle of
a comodule $M\in\mathcal{M}^{\Gamma}$ is denoted by soc$(M)$. If $M=\Gamma$
then soc$(M)$=corad($\Gamma),$ the coradical. The \textit{coefficient space}
of $M$ in $\Gamma$ is denoted by $\mathrm{cf}(M).$

$M$ is said to be \textit{quasifinite} if $\operatorname*{Hom}_{-\Gamma}(F,M)$
is finite-dimensional for all $F\in\mathcal{M}_{f}^{\Gamma}.$ Equivalently,
the simple summands of soc$(M)$ have finite (but perhaps unbounded)
multiplicities [Ch].

An \textit{injective copresentation} of $M$ is an exact sequence $0\rightarrow
M\rightarrow I_{0}\rightarrow I_{1}$, where $I_{0}$ and $I_{1}$ are
injectives. The comodule $M$ is said to be \textit{quasifinitely copresented}
if $I_{0}$ and $I_{1}$ can be chosen to be quasifinite.

Assume for the moment that $\Gamma$ is a finite-dimensional, so that
$\mathcal{M}_{f}^{\Gamma}\approx\operatorname{mod}R.$ Let $\operatorname*{tr}$
denote the usual transpose on $\operatorname{mod}R$, as in the representation
theory of finite-dimensional algebras (defined using projective resolutions,
see [ARS, \S IV]). It is easy to see from duality that for finite-dimensional
comodules that \textrm{D}$\operatorname*{Tr}$\textrm{D}$=\operatorname*{tr}$,
so we can say \textrm{D}$\operatorname*{tr}=\operatorname*{Tr}$\textrm{D} and
$\operatorname*{tr}$\textrm{D}$=$\textrm{D}$\operatorname*{Tr}$.

\section{Duality for quasifinite injective comodules}

We define the contravariant functor $\ast$ as in [CKQ] as $h_{-\Gamma
}(-,\Gamma):\mathcal{M}_{q}^{\Gamma}\rightarrow\mathcal{M}^{\Gamma}$ (also the
version on the left as well)

With $R=\mathrm{D}\Gamma,$ we have the right and left hit actions of $R$ on
$\Gamma$ (more generally on any $\Gamma,\Gamma$-bicomodule), usually denoted
by the symbols $\leftharpoonup,\rightharpoonup$ as in e.g. [Mo, Sw]. Here we
will omit these symbols and simply use juxtapostion, e.g., $e\Gamma
=e\rightharpoonup\Gamma$, $e\in R$. Notice that $\Gamma e$ is an injective
right comodule.

\begin{theorem}
[CKQ]$\ast$ restricts to a duality on the category of quasifinite injective
comodules $\mathcal{I}^{\Gamma}.$
\end{theorem}

We first collect some facts concerning duality and injectives.

\begin{lemma}
Let $e=e^{2}\in R$. Then:\newline (a) $\mathrm{D}(\Gamma e)\cong eR$ as right
$R$-modules.\newline (b) If $\Gamma\mathrm{\ }$is of finite dimension, then
\textrm{D}$(eR)\cong\Gamma e$ as left $R$-modules\newline (c)
$\operatorname*{Hom}_{-\Gamma}(\Lambda,\Gamma e)=\operatorname*{Hom}_{-\Gamma
}(\Lambda,\Lambda e)$ for every subcoalgebra $\Lambda\subset\Gamma\newline
$(d) $(\Gamma e)^{\ast}\cong e\Gamma$ as left $\Gamma$-comodules.
\end{lemma}

\begin{proof}
The statement in (a) is essentially the definition of the right hit action as
dual to the left multiplication by $R.$ We leave the details to the reader.
The proof of (b) follows from (a) and the duality D for finite dimensional (co)modules.

For (c), observe that the image of any comodule map on the left-hand term has
its image in $\Delta^{-1}(\Gamma e\otimes\Lambda)\cong\Gamma e\square\Lambda.$
The result follows immediately from the additivity of the cotensor product, as
$\Gamma e$ is a summand of $\Gamma.$

Lastly we prove (d). Let $\Lambda$ denote a finite-dimensional subcoalgebra of
$\Gamma.$ Then we have
\begin{align}
Hom_{-\Gamma}(\Lambda,\Gamma e)  &  =\operatorname*{Hom}\nolimits_{-\Gamma
}(\Lambda,\Lambda e)\nonumber\\
&  \cong\operatorname*{Hom}\nolimits_{R}(\mathrm{D}(\Lambda e),\mathrm{D}%
\Lambda)\nonumber\\
&  \cong\operatorname*{Hom}\nolimits_{R}(e(\mathrm{D}\Lambda),\mathrm{D}%
\Lambda)\nonumber\\
&  \cong(\mathrm{D}\Lambda)e
\end{align}
by applying the part (c), the duality $\mathrm{D,}$ and part (a). Note too
that we consider $e$ to act on \textrm{D}$\Lambda$ via the ring map
\textrm{D}$\Gamma\rightarrow\mathrm{D}\Lambda$ given by restriction. Now we
obtain
\begin{align*}
(\Gamma e)^{\ast}  &  =\lim\limits_{\rightarrow}\mathrm{D}Hom_{-\Gamma
}(\Lambda,\Gamma e)\\
&  \cong\lim\limits_{\rightarrow}\mathrm{D}((D\Lambda)e)\\
&  \cong\lim\limits_{\rightarrow}e\Lambda\\
&  \cong e\Gamma
\end{align*}
where the direct limits are over the finite-dimensional subcoalgebras
$\Lambda$, using the definition of $\ast$, (1) above, and part (b). This
completes the proof of the Lemma.
\end{proof}

\bigskip

\begin{proof}
(of the Theorem) Let $I$ be a quasifinite injective comodule. Then I is the
direct sum of indecomposable injectives, all of the form $\Gamma e$, where is
a primitive idempotent in $R$. Part (d) of the Lemma and its right-handed
counterpart yield $I^{\ast\ast}\cong I,$ an isomorphism that is easily seen to
be natural.
\end{proof}

\section{Direct Limits}

\subsection{The Transpose}

Let $M\in\mathcal{M}_{qc}^{\Gamma}$. The transpose $\operatorname*{Tr}M\in$
$^{\Gamma}\mathcal{M}$ is defined in [CKQ] to be by $0\rightarrow
\operatorname*{Tr}M\rightarrow I_{0}^{\ast}\rightarrow I_{1}^{\ast}$ where
$0\rightarrow M\rightarrow I_{0}^{{}}\rightarrow I_{1}^{{}}$is a minimal
quasifinite injective copresentation of $M$. Actually, it is easy to see that
that we can use any quasifinite injective copresentation of $M$ to define
$\operatorname*{Tr}M$, i.e., we may omit the minimality requirement from the
definition. \ For a finite-dimensional subcoalgebra $\Lambda\subset\Gamma,$
$\operatorname*{Tr}_{\Lambda}M$ denotes the transpose of $M\in\mathcal{M}%
_{qc}^{\Lambda}$

\begin{lemma}
\ Let $M$ $\in\mathcal{M}_{f}^{\Gamma}$ be a finite-dimensional quasifinitely
copresented comodule. Then for every finite-dimensional subcoalgebra
$\Lambda\subset\Gamma$ containing $\operatorname*{cf}(M)$ \
\[
\operatorname*{Tr}\nolimits_{_{\Lambda}}M\cong\Lambda\square\operatorname*{Tr}%
M
\]
\end{lemma}

\begin{proof}
Let $\Lambda$ be as in the statement. Let
\[
0\rightarrow M\rightarrow I_{0}\Box\Lambda\rightarrow I_{1}\Box\Lambda
\]
be an injective copresentation of $M$. We obtain an injective copresentation
\[
0\rightarrow M\rightarrow I_{0}\Box\Lambda\rightarrow I_{1}\Box\Lambda
\]
for $M$ in $\mathcal{M}^{\Lambda}$. \ The defining copresentation for the
transpose in $\mathcal{M}^{\Lambda}$ using ``$\ast$ in $\mathcal{M}^{\Lambda}
$'' is
\[
0\rightarrow\operatorname*{Tr}\nolimits_{_{\Lambda}}M\rightarrow h_{-\Lambda
}(I_{1}\Box\Lambda,\Lambda)\rightarrow h_{-\Lambda}(I_{0}\Box\Lambda,\Lambda).
\]
On the other hand, the defining copresentation
\[
0\rightarrow\operatorname*{Tr}M\rightarrow h_{-\Gamma}(I_{1},\Gamma
)\rightarrow h_{-\Gamma}(I_{0},\Gamma)
\]
for $M$ in $\mathcal{M}^{\Gamma}$ can be cotensored with $\Lambda$, yielding
\[
0\rightarrow\Lambda\square\operatorname*{Tr}M\rightarrow\Lambda\square
h_{-\Gamma}(I_{1},\Lambda^{\prime})\rightarrow\Lambda\square h_{-\Gamma}%
(I_{0},\Lambda)
\]
In view of [Tak, 1.14(b)], which says that $\Lambda\Box h_{-\Gamma}%
(M,\Gamma)\cong h_{-\Lambda}(M,\Lambda),$ the result is established.
\end{proof}

Assume that $\Lambda$ is a finite-dimensional coalgebra. We see next that the
transpose commutes with certain direct limits: Let $\operatorname*{tr}%
_{D\Lambda}$ and $\operatorname*{Tr}_{\Lambda}$ denote the transposes in the
module category $_{D\Lambda}\mathcal{M}$ and the comodule category
$\mathcal{M}^{\Lambda}$, respectively. Since cotensoring commutes with direct
limits and is left exact, it follows from the preceding lemma that

\begin{proposition}
Suppose that $M\in\mathcal{M}_{qc}^{\Gamma}$ finite-dimensional,
indecomposable and not injective. $\ $\ Then
\[
\lim\limits_{\overrightarrow{\mathcal{L}}}\operatorname*{Tr}\nolimits_{\Lambda
}M=\operatorname*{Tr}M
\]
\ for the direct system $\mathcal{L}$ \ of finite-dimensional
subcoalgebras$\ $of $\ \Gamma$.
\end{proposition}

\begin{remark}
The preceding result gives a way of constructing indecomposables which are
unions of chains of indecomposable subcomodules. See Example 3 at the end of
this paper for an example of such an infinite-dimensional indecomposable.
\end{remark}

We show next that the almost split sequences (e.g. those whose existence is
provided by [CKQ]) are direct limits of almost split sequences over
finite-dimensional subcoalgebras.

\subsection{Almost Split Sequences as Direct Limits}

\begin{proposition}
\ Suppose that $C\in\mathcal{M}_{qc}^{\Gamma}$ is finite-dimensional. \ \ Let
\[
d:0\rightarrow A\rightarrow B\rightarrow C\rightarrow0
\]
be an almost split sequence in $\mathcal{M}^{\Gamma}$. \ Then $d$ is the
direct limit of almost split sequences $(d\Box\Lambda)_{\Lambda\in\mathcal{L}%
}$ for some direct system $\mathcal{L}$ \ of finite-dimensional subcoalgebras
of $\Gamma$.
\end{proposition}

\begin{proof}
Let $\Lambda$ be any finite-dimensional coalgebra of $\Gamma$ containing both
the coefficient space of $C$, and a finite-dimensional subspace of $B$ mapping
onto $C.$ \ It is clear that the sequence $d\Box\Lambda:$ $B\Box
\Lambda\rightarrow C\Box\Lambda=C\rightarrow0$ is exact (onto $C$) for any
finite-dimensional subcoalgebra $\Lambda\subseteq\Gamma$, provided
$\Lambda\supseteq$ $\Lambda^{\prime}$. \ \ It is trivial to check that
$B\Box\Lambda\rightarrow C\rightarrow0$ is right almost split in
$\mathcal{M}^{\Lambda}.$

The direct system $\mathcal{L}$ of subcoalgebras is taken to be the
finite-dimensional subcoalgebras $\Lambda$ as just described, ordered by
inclusion. The Lemma ensures that $\ A\Box\Lambda=$\textrm{$\operatorname*{Tr}%
$}$_{\Lambda}D(C)$, and \textrm{$\operatorname*{Tr}$}$_{\Lambda}D(C)$ is
indecomposable by the finite-dimensional theory.\ This is the left-hand term
of the right almost split sequence $d\Box\Lambda$, so by the standard result
of Auslander (see [CKQ], 4.3 Proof) $\delta\Box\Lambda$ is in fact an almost
split sequence in $\mathcal{M}^{\Lambda}$. \ The maps between the sequences is
the obvious one using inclusions.
\end{proof}

\section{Finitely Almost Split Sequences}

\subsection{\textbf{The Mittag-Leffler Condition}}

For abelian groups, the direct limit (or filtered colimit) is an exact functor
but (inverse) limits are left exact but not always right exact. The
Mittag-Leffler condition [Gro] (see also [Ha, p.119]) guarantees that inverse
limits of certain inverse systems of exact sequences are exact. \ It has the
following special case:

\begin{proposition}
Let
\[
d_{i}:0\rightarrow M_{i}\rightarrow E_{i}\rightarrow N_{i}\rightarrow0
\]
be an inverse system of exact sequences of vector spaces where the index set
has a countable cofinal subset. If each $M_{i}$ is finite-dimensional, then
\[
\lim\limits_{\leftarrow}d_{i}=d:0\rightarrow M\rightarrow E\rightarrow
N\rightarrow0
\]
is exact.
\end{proposition}

The reader may consult the references just mentioned, or [Je, Proposition 2.3]
for a simple proof.

\subsection{Prorational Modules}

Define $\mathrm{D}:$ $^{\Gamma}\mathcal{M}\rightarrow\mathcal{D}^{\Gamma}$ to
be the functor $\mathrm{D}$ with range being the full subcategory of
$\mathcal{M}_{R}$ having objects $\mathrm{D}M$, $M$ $\in$ $^{\Gamma
}\mathcal{M}$. \ We let $\mathcal{D}_{q}^{\Gamma}$ denote the full subcategory
of duals of quasifinite comodules. Note that the objects in $\mathcal{D}%
^{\Gamma}$ are precisely the inverse limits of finite-dimensional rational
right $R$-modules$.$ Accordingly the finite-dimensional objects in
$\mathcal{D}^{\Gamma}$ are precisely the finite-dimensional rational right
$R$-modules$.$ We refer to the objects of $\mathcal{D}^{\Gamma}$ as
\textit{prorational} right $R$- modules. \ We shall freely use the variant of
$\mathrm{D}$ (denoted by the same symbol) with opposite categories, e.g.
$\mathrm{D}:$ $\mathcal{M}^{\Gamma}\rightarrow$ $^{\Gamma}\mathcal{D}$.

\begin{definition}
We say that $\mathrm{D}M$ $\in$ $\mathcal{D}^{\Gamma}$ is
\textbf{coquasifinite }if $\ \operatorname*{Hom}_{R}(\mathrm{D}M,F)$ is
finite-dimensional for all finite-dimensional $F\in$ $\mathcal{D}^{\Gamma}$.
\end{definition}

The dual notion ``quasifinite'' was introduced in [Tak]. The objects of
$\mathcal{D}_{q}^{\Gamma}$ are coquasifinite prorational right $R$-modules.

We show below that $\mathcal{D}^{\Gamma}$ (resp. $\mathcal{D}_{q}^{\Gamma}) $
is the dual category of the category of (resp. quasifinite) comodules. \ 

If $\rho:M\rightarrow M\otimes\Gamma$ is the structure map of $M$ $\in$
$\mathcal{M}^{\Gamma}$, then by restricting to $\mathrm{D}M\otimes
R\subset\mathrm{D}(M\otimes R)$ and abusing notation, we obtain the map
\[
\mathrm{D}\rho:\mathrm{D}M\otimes R\rightarrow\mathrm{D}M
\]
It is straightforward to check that this coincides with the rational right $R
$-module structure on $\mathrm{D}M$ which arises from the left $\Gamma
-$comodule structure on \textrm{D}$\mathrm{D}M$, which in turn comes from
$\rho.$

Note that the objects in $\mathcal{D}^{\Gamma}$ are precisely the inverse
limits of finite-dimensional rational right $R$-modules$.$ Accordingly the
finite-dimensional objects in $\mathcal{D}^{\Gamma}$ are precisely the
finite-dimensional rational right $R$-modules (i.e. left comodules). We refer
to the objects of $\mathcal{D}^{\Gamma}$ as \textit{prorational modules.}

In the following lemma, a finite-dimensional cotensor product of comodules is
seen to be dual to the tensor product.

\begin{lemma}
Let $M\in\mathcal{M}^{\Gamma}$ and $N\in$ $^{\Gamma}\mathcal{M}^{{}}$ and
assume that $M\Box N$ is finite-dimensionsal. Then $\mathrm{D}(M\Box N)$ is
isomorphic to $\mathrm{D}M\otimes_{R}\mathrm{D}N$.
\end{lemma}

\begin{proof}
Let $\rho:M\rightarrow M\otimes\Gamma$ and $\lambda:N\longrightarrow
\Gamma\otimes N$ be the structure maps of $M$ and $N$ respectively. \ Then as
is noted above, $\mathrm{D}M$ is a left comodule and is a (rational) right
$R$-module. \ Similarly $\mathrm{D}N$ is a left $R$-module.

The cotensor $M\Box N$ is defined by the usual equalizer
\[
M\Box N\rightarrow M\otimes N\overset{\rho\otimes1}{\underset{1\otimes\lambda
}{\rightrightarrows}}M\otimes\Gamma\otimes N.
\]
\ Dualizing, we have the coequalizer
\[
\mathrm{D}(M\otimes\Gamma\otimes N)\overset{\mathrm{D}(\rho\otimes
1)}{\underset{\mathrm{D}(1\otimes\lambda)}{\rightrightarrows}}\mathrm{D}%
(M\otimes N)\overset{p}{\longrightarrow}\mathrm{D}(M\Box_{\Gamma
}N)\longrightarrow0.
\]

By hypothesis $M\Box_{\Gamma}N$ is finite-dimensional, so the density of
$\mathrm{D}(M)\otimes\mathrm{D}(N)$ in $\mathrm{D}(M\otimes N)$ implies that
the restriction $\mathrm{D}M\otimes\mathrm{D}N\longrightarrow\mathrm{D}(M\Box
N)$ is onto. \ The kernel of this map is $\ker p\cap(\mathrm{D}M\otimes
\mathrm{D}N).$ Thus we have the coequalizer
\[
\mathrm{D}M\otimes\mathrm{D}\Gamma\otimes DN\overset{\mathrm{D}\rho\otimes
1}{\underset{1\otimes\mathrm{D}\lambda}{\rightrightarrows}}\mathrm{D}M\otimes
DN\overset{p}{\longrightarrow}\mathrm{D}(M\Box_{\Gamma}N).
\]
This finishes the proof of the lemma.
\end{proof}

Let $\Lambda\ $denote a finite-dimensional subcoalgebra of $\Gamma$. The
finite-dimensional algebra $\mathrm{D}\Lambda$ is isomorphic to $\mathrm{D}%
\Gamma/\Lambda^{\perp}$, where $\Lambda^{\perp}$ is the ideal of functionals
in $R=\mathrm{D}\Gamma$ vanishing on $\Lambda$. \ The following is now
immediate. Let $\mathcal{L}$ denote the direct system of finite-dimensional
subcoalgebras of $\Gamma.$ Lemma 2 immediately yields

\begin{proposition}
Let $M_{\Lambda}$ denote $M\Box\Lambda$, $M$ $\in$ $\mathcal{M}_{q}^{\Gamma}%
$.\ \ Then \newline (a) $\mathrm{D}(M_{\Lambda})\cong\mathrm{D}M\otimes
_{R}\mathrm{D}\Lambda\cong\mathrm{D}M/M\Lambda^{\perp}.$ \newline (b)
$\mathrm{D}M$ $\cong\lim\limits_{\underset{\mathcal{L}}{\leftarrow}}%
\mathrm{D}M_{\Lambda}$.
\end{proposition}

\begin{proposition}
(a) $\mathrm{D}$: $^{\Gamma}\mathcal{M}\rightarrow\mathcal{D}^{\Gamma}$ is a
duality \newline (b) \textrm{D} restricts to a duality $^{\Gamma}%
\mathcal{M}_{q}\rightarrow\mathcal{D}_{q}^{\Gamma}$\newline (c) $\mathrm{D}M$
is coquasifinite for all $M$ $\in$ $^{\Gamma}\mathcal{M}_{q}.$
\end{proposition}

\begin{proof}
A functor \textrm{D}$^{\prime}$:$\mathcal{D}^{\Gamma}\rightarrow$ $^{\Gamma
}\mathcal{M}$ can be defined by setting \textrm{D}$^{\prime}$\textrm{D}$M$ =
$M$, giving a correspondence on objects. Suppose $f\in\operatorname*{Hom}%
_{R}($\textrm{D}$N,$\textrm{D}$M).$ \ We put \textrm{D}$^{\prime}%
f=\lim\limits_{\rightarrow}$\textrm{D}$f_{i}$ where $f_{i}=f\otimes R_{i}.$
\ One can check that the functor \textrm{D}$^{\prime}$ is an inverse duality
for \textrm{D}$. $ We are done with (a). Part (c) follows from the duality,
and part (b) is clear.
\end{proof}

\subsection{Finitely Almost Split Sequences}

\begin{definition}
We say that a nonsplit exact sequence of objects $0\rightarrow A\rightarrow
B\rightarrow C\rightarrow0$ in $_{R}\mathcal{M}$ is \textbf{right finitely
almost split} if

1. $C$ is indecomposable

2. every morphism $X\rightarrow C$ with $X$ $\in\mathcal{M}_{f}^{\Gamma}$ ,
which is not a split epimorphism, lifts to $B$. \ 

Dually, the sequence is said to be \textbf{left finitely almost split} if $A$
is indecomposable if every morphism $A\rightarrow X$ with $X$ $\in
\mathcal{M}_{f}^{\Gamma},$ which is not a split monomorphism, extends to $B$.
\end{definition}

These definitions are a coalgebraic version of the definition of finitely
almost split sequences given for ``local nests of quivers'' given in [GM],
where they are called ``special sequences''. \ Further extending their work we have

\begin{theorem}
Let $\Gamma$ be a countable-dimensional coalgebra. If $C\in\mathcal{M}%
_{f}^{\Gamma}$ is nonprojective and indecomposable, then there exists a right
finitely almost split sequence
\[
d:0\mathcal{\rightarrow}A\mathcal{\rightarrow}B\mathcal{\rightarrow
}C\mathcal{\rightarrow}0
\]
in $\mathcal{M}_{q}^{\Gamma}$. \ If $\mathrm{D}C$ is quasifinitely
copresented, then $d$ can be chosen to be an almost split sequence in
$\mathcal{M}_{q}^{\Gamma}$.
\end{theorem}

\begin{proof}
We can write $\Gamma=\cup\Gamma_{i}$ as the ascending union of a chain of
finite-dimensional subcoalgebras (not necessarily the coradical filtration)
$\Gamma_{i}$ where we may assume that $C\Box\Gamma_{0}=C.$

It is easy to check that $C$ is not projective in $\mathcal{M}^{\Lambda_{i}} $
for all $i$. For, let $p:B\rightarrow C$ be a surjection in $\mathcal{M}%
^{\Gamma};$ then $B\square\Lambda_{i}\rightarrow C\square\Lambda_{i}=C$ is a
surjection in $\mathcal{M}^{\Lambda_{i}}$. Were $C$ projective in
$\mathcal{M}^{\Lambda_{i}}$, then a splitting map $C\rightarrow B\square
\Lambda_{i}$ in $\mathcal{M}^{\Lambda_{i}}$ would also split $p.$

Let $\operatorname*{tr}_{i}$ and $\operatorname*{Tr}_{i}$ denote the
transposes in the module category $_{\mathrm{D}\Lambda_{i}}\mathcal{M}$ and in
$\mathcal{M}^{\Lambda_{i}}$, respectively, (as mentioned in the introduction).
By results for finite-dimensional algebras, there exists an almost split
sequence
\[
d_{i}:0\rightarrow A_{i}\rightarrow B_{i}\rightarrow C\rightarrow0
\]
in $\mathcal{M}_{f}^{\Gamma_{i}}$ with $A_{i}$=\textrm{D}$\operatorname*{tr}%
_{i}C$=$\operatorname*{Tr}_{i}$\textrm{D}$C,$ $i\in\mathbb{N}$. We obtain
exact sequences indexed by $\mathbb{N}$, along with maps $d_{i}\rightarrow
d_{i+1}$ obtained by assigning the identity map on $C$, and then a little
diagram chasing using the fact that $d_{i+1}$ is almost split. Furthermore,
the maps $A_{i}\rightarrow A_{i+1}$ and $B_{i}\rightarrow B_{i+1}$are
monomorphisms. This can be seen as follows. Suppose that $A_{i}\rightarrow
A_{i+1}$ has kernel $K_{i}\neq0$. Then we obtain an exact sequence
\[
d_{i}^{\prime}:0\rightarrow A_{i}^{\prime}\rightarrow B_{i}^{\prime
}\rightarrow C\rightarrow0
\]
where $A_{i}^{\prime}=A_{i}/K$ and $B_{i}^{\prime}=B_{i}/K$, with epimorphisms
$d_{i}\rightarrow d_{i}^{\prime}$ (identity map on $C$). It follows easily
from the fact that $d_{i}$ is almost split that $d_{i}^{\prime}$ is split. But
this immediately implies that $d_{i}$ splits. Thus we obtain a direct system
of exact sequences $(d_{i})_{i\in\mathbb{N}}$ with monomorphic connecting maps
$d_{i}\rightarrow d_{i+1}$, each being the identity on $C$.

The direct limit is an exact sequence
\[
d:0\mathcal{\rightarrow}A\mathcal{\rightarrow}B\mathcal{\rightarrow
}C\mathcal{\rightarrow}0.
\]
since the direct limit is an exact functor.

We show that $d$ is right finitely almost split. The sequence is not split,
for otherwise a splitting map $C\rightarrow B$ would have its image in $B_{i}
$ for some $i$ (since $C$ is finite-dimensional). This would then be a
splitting of $d_{i}$, a contradiction. Let $X$ $\in\mathcal{M}_{f}^{\Gamma} $,
so that $X\in$ $\mathcal{M}^{\Gamma_{i}}$ for some $i$. Then any morphism $X$
$\rightarrow C$, which is not a split epimorphism, lifts to $B_{i}$. Composing
with the inclusion map $B_{i}\rightarrow B$, we get the required lifting. This
shows that $d$ is right finitely almost split.

Now assume as in the statement that \textrm{D}$C$ is quasifinitely
copresented. To show that $d$ is right almost split in $\mathcal{M}%
_{q}^{\Gamma}$, it suffices to show that the sequence $\operatorname*{Hom}%
(X,d)$
\[
0\rightarrow\operatorname*{Hom}(X,A)\rightarrow\operatorname*{Hom}%
(X,B)\rightarrow\operatorname*{Hom}(X,C)\rightarrow0
\]
is exact for all $X\in$ $\mathcal{M}_{q}^{\Gamma}$. Let $X_{i}=X\Box\Gamma
_{i}$, which is a finite-dimensional comodule (since $X$ is quasifinite). Note
that $\operatorname*{Hom}(X_{i},d)$ is an inverse system of short exact
sequences. By the Mittag-Leffler condition (Proposition 3), it suffices to
show $\operatorname*{Hom}(X_{i},A)$ is finite-dimensional for all $i$. Thus it
suffices to know that $A$ is quasifinite. First, we have $A=\lim
\limits_{\rightarrow}A_{i}$=$\lim\limits_{\rightarrow}\operatorname*{Tr}_{i}%
$\textrm{D}$C$ =$\operatorname*{Tr}$\textrm{D}$C$ (3.1 Lemma 2). Secondly, the
rudimentary property of $\operatorname*{Tr}$ [CKQ, Lemma 3.2] states in fact
that $\operatorname*{Tr}$\textrm{D}$C$ is quasifinite (because \textrm{D}$C$
is quasifinitely copresented).

By duality, \textrm{D}$C$ in indecomposable and noninjective. By [CKQ, 3.2]
and 3.1 Lemma 2, $\ A=\operatorname*{Tr}$\textrm{D}$C$ is also indecomposable
noninjective. \ By a result of Auslander (see [CKQ, 4.3]), $d$ is an almost
split sequence.

This completes the proof of the Theorem.
\end{proof}

\begin{remark}
The second statement of the theorem holds ([CKQ, Corollary 4.3]) without the
assumption on dimension, and the almost sequences obtained are almost split in
all of $\mathcal{M}^{\Gamma}$
\end{remark}

The argument of the theorem dualizes.

\begin{theorem}
Let $\Gamma$ be a countable-dimensional coalgebra. Let $A\in\mathcal{M}%
_{f}^{\Gamma}$. If $A$ is noninjective and indecomposable, then there exists a
left finitely almost split sequence
\[
d:0\mathcal{\rightarrow}A\mathcal{\rightarrow}B\mathcal{\rightarrow
}C\mathcal{\rightarrow}0
\]
in $\mathcal{D}_{q}^{\Gamma}$. If $A$ is quasifinitely copresented, then $d$
can chosen to be an almost split sequence in $\mathcal{D}_{q}^{\Gamma}$.
\end{theorem}

\begin{proof}
Since the proof is dual to the one above we only give sample details. Assume
as in the statement that $A$ is quasifinitely copresented and we have a left
finitely almost split sequence
\[
d:0\mathcal{\rightarrow}A\mathcal{\rightarrow}B\mathcal{\rightarrow
}C\mathcal{\rightarrow}0
\]
in $\mathcal{D}_{q}^{\Gamma}$. Here $d$ is obtained as the inverse limit of
sequences $d_{i}:0\mathcal{\rightarrow}A\mathcal{\rightarrow}B_{i}%
\mathcal{\rightarrow}C_{i}\mathcal{\rightarrow}0$ , that are almost split in
$\mathcal{D}_{q}^{\Gamma_{i}}$. To show that $d$ is left almost split in
$\mathcal{M}_{q}^{\Gamma}$, it suffices to show that the sequence
$\operatorname*{Hom}(d,Y)$ is exact for all $Y\in$ $\mathcal{D}_{q}^{\Gamma}$.
We have $Y=\lim\limits_{\longleftarrow}$ $Y_{i}$ where the $Y_{i}$ are finite
dimensional rational $R$-modules. According to the Mittag-Leffler condition,
it suffices to show $\operatorname*{Hom}(A,Y_{i})$ is finite-dimensional for
all $i $. Thus we want to know that $A$ is coquasifinite. First, we have
\[
C=\lim\limits_{\longleftarrow}C_{i}=\lim\limits_{\longleftarrow}%
\mathrm{D}\operatorname*{Tr}\nolimits_{i}A=\mathrm{D}\lim\limits_{\rightarrow
}CA=\mathrm{D}\operatorname*{Tr}A
\]
Secondly [CKQ, Lemma 3.2] states that $\operatorname*{Tr}A$ is quasifinite.
Therefore $C=\mathrm{D}\operatorname*{Tr}A$ is coquasifinite.

By [CKQ, 3.2] and 3.1 Lemma 2, $\operatorname*{Tr}A$ is an indecomposable
noninjective comodule. \ Hence $C=\mathrm{D}\operatorname*{Tr}A$ is
indecomposable and nonprojective in $\mathcal{D}_{q}^{\Gamma}$.
\end{proof}

By dualizing [CKQ], Corollary 4.3(b), we obtain a result extending the second
statement of the Thereom above. We leave the details to the reader.

\begin{theorem}
Let $A\in\mathcal{M}_{f}^{\Gamma}$. If $A$ is quasifinitely copresented,
indecomposable and not injective, then exists an almost split sequence
$0\rightarrow A\rightarrow B\rightarrow C\rightarrow0$ in $\mathcal{D}%
^{\Gamma}$ with $C=\mathrm{D}\operatorname*{Tr}A$.
\end{theorem}

\begin{remark}
With $A$ in the hypothesis, \textrm{D}$A$ is coquasifinitely presented. This
means that there exists a projective resolution \textrm{D}$I_{1}%
\rightarrow\mathrm{D}I_{0}\rightarrow\mathrm{D}A\rightarrow0$ in
$\mathcal{D}^{\Gamma},$ with the \textrm{D}$I_{i}$ coquasifinite. \ If $I$ is
a finitely cogenerated injective comodule (see [Ch]), then $\mathrm{D}I$ is a
projective as an R-module. \ But if $I$ is only assumed to be quasifinite
injective, we do not know this to be the case.
\end{remark}

\begin{remark}
The construction of almost split sequences as limits gives a way of showing
that the category of finite-dimensional comodules does not have almost split
sequences. \ Let $C$ be a finite-dimensional comodule such that
$\operatorname*{Tr}\mathrm{D}C$ is infinite-dimensional. \ Then there is no
almost split sequence in $\mathcal{M}_{f}^{\Gamma}$ ending in $C$. \ Dually,
if $\operatorname*{Tr}A$ is infinite-dimensional, then there is no almost
split sequence starting with $A$ in $\mathcal{M}_{fd}^{\Gamma}.$ An example of
this type is given below. On the other hand, if $\mathcal{M}_{f}^{\Gamma}$ has
almost split sequences (e.g. if $\Gamma$ is right semiperfect, see [CKQ]),
then the Auslander-Reitan quiver exists. Examples of AR-quivers for
$\mathcal{M}_{f}^{\Gamma}$ for various path coalgebras are given in [NS].
\end{remark}

\section{Examples}

The coalgebras in the next two examples are neither right nor left
semiperfect, and finite-dimensional indecomposables have finite-dimensional
transposes. The almost split sequences are almost split in the category of
finite-dimensional comodules. \ In these examples, the almost split sequence
obtained is the same as the one obtained over a finite-dimensional
subcoalgebra (i.e.\ over a finite subquiver).

In contrast, the third example shows that a simple comodule can have an
infinite-dimensional transpose.\ So the category of finite-dimensional
comodules does not have almost split sequences.

\bigskip

1. Let $Q$ be the quiver of type $\mathbb{A}_{\infty}^{\infty}$ with vertices
labeled by the integers and arrows $a_{i}:i\rightarrow i+1,$ $i\in\mathbb{Z}$
. \ We write $S(i)$ for the simple in $\mathcal{M}^{kQ}$ corresponding to each
vertex $i$ and denote its injective hull by $I(i)$

By the theory of Nakayama algebras and Dynkin quivers (see [ARS] and [Ga]),
the isomorphism classes of finite-dimensional indecomposable comodules are
given by the representations $V_{ij}=$ $(V,f)_{ij}$ ( for all $i\leq j$)
defined by $V_{t}=$k for $i\leq t\leq j$, zero otherwise; the linear maps are
$f_{t}:V_{t}\rightarrow V_{t+1}$ given by $f_{t}=1$ for $i\leq t<j$ and zero
otherwise. \ We compute injective envelope of $V_{ij}$ to be
\[
I(V_{ij})=I(j)
\]
and we find that $I(V_{ij})/V_{ij}\cong I(i-1)$. \ Thus
\[
0\rightarrow\operatorname*{Tr}(V_{ij})\rightarrow I(i-1)^{\ast}\rightarrow
I(j-1)^{\ast}%
\]
is a copresentation yielding \textrm{D}$\operatorname*{Tr}(V_{ij}%
)=V_{i-1,j-1}$. \ The almost split sequences are
\[
0\rightarrow V_{i,j}\rightarrow V_{i-1,j}\oplus V_{i,j-1}\rightarrow
V_{i-1,j-1}\rightarrow0
\]
with irreducible maps $V_{i,j}\rightarrow V_{i,j+1}$ and $V_{i-1,j}\rightarrow
V_{i,j}$ being the obvious monomorphism into the first summand and epimorphism
onto the second summand. The map on the right is given by natural epimorphism
and monomorphism with alternate signs.

\bigskip

2. \ Similarly, let $Q$ denote the quiver with one vertex and one loop. \ Then
there is a unique finite-dimensional indecomposable right comodule $V_{n}$ of
dimension $n\geq0$. \ It is straightforward to see that \textrm{D}%
$\operatorname*{Tr}V_{n}=V_{n}$ and that the almost split sequences are given
by
\[
d_{n}:0\rightarrow V_{n}\rightarrow V_{n-1}\oplus V_{n+1}\rightarrow
V_{n}\rightarrow0
\]
(just as in [ARS, p. 141]). Following the ideas in this article. one might
hope to take the limit of these sequences to get and almost split sequence for
the infinite-dimensional indecomposable $V=\lim\limits_{\rightarrow}V_{n}\cong
kQ$. Unfortunately, the sequence so obtained is split.

\bigskip

3. Let $Q$ be the quiver of type $\mathbb{D}_{\infty}$ with vertices labeled
by positive integers and two special vertices $0,0\prime$. The arrows are
defined to be $a_{i}:i\rightarrow i+1,$ $i\in\mathbb{Z}^{+}$, and
$a_{0}:0\rightarrow1,$ $a_{0^{\prime}}:0^{\prime}\rightarrow1$. \ 

We compute the transpose of the simple right noninjective comodule $S(1).$ \ 

$I(1)$ is given by the representation $V=$ $(V,f)$ defined by $V_{t}%
=\mathrm{k}$ for $t=1,0,0^{\prime}$ and zero otherwise, with linear maps are
$f_{t}:V_{t+1}\rightarrow V_{t}$ given by $f_{t}=1$ for $t=0,0^{\prime}$. Thus
$I(1)/S(1)$ is $S(0)\oplus S(0)^{\prime},$ the direct sum of two simple injectives.

Next observe that $I(S(0)\oplus S(0)^{\prime})^{\ast}$ = $I(0)^{\ast}\oplus
I(0^{\prime})^{\ast}$ is given by the representation $V=$ $(V,f)$ with by
$V_{t}=\mathrm{k}$ for $t=0,0^{\prime}$, $V_{t}=\mathrm{k}^{2}$ for $t>0$, and
zero otherwise. \ We leave the maps to the reader.

$I(1)^{\ast}$is given by the representation $V=$ $(V,f)$ with by $V_{t}=0$ for
$t=0,0^{\prime}$, $V_{t}=\mathrm{k}$ for $t>0$, and zero otherwise.

Finally, $0\rightarrow\operatorname*{Tr}(S(1))\rightarrow I(S(0)\oplus
S(0)^{\prime})^{\ast}\rightarrow I(1)^{\ast}$ is a copresentation yielding
$\operatorname*{Tr}(1)=V$ where the representation $V=$ $(V,f)$ of $Q^{op}$ is
defined by $V_{t}=\mathrm{k}$ for all vertices $t$, with linear maps are
$f_{t}:V_{t+1}\rightarrow V_{t}$ given by $f_{t}=1$ for all $t$. $\ $\ We know
that $V$ is an infinite-dimensional indecomposable left comodule; thus
$\mathcal{M}_{f}^{\mathrm{k}Q}$ does not have an almost split sequence
\textit{starting} at $S(1)$. \ The almost split sequence in $\mathcal{D}%
_{q}^{kQ}$ given by the Theorem 4 is of the form
\[
0\rightarrow S(1)\rightarrow B\rightarrow\mathrm{D}V\rightarrow0
\]
$.$

\end{document}